\documentstyle[txmac,a4,
amssymb,%
case,%
twoside,%
nocaphead,%
epsf,%
varthm,%
rotate,%
myrot,%
mypic,times,mathptm]{article}

\advance\oddsidemargin by -1.9cm
\advance\evensidemargin by -1.9cm
\advance\textwidth by 3.8cm

\def\mynewtheo#1#2{%
\newtheorem{@#1}{#2}[section]%
\newenvironment{#1}{\begin{@#1}\rm}{\end{@#1}}}

\mynewtheo{lemma}{Lemma}
\mynewtheo{exer}{Exercise}
\mynewtheo{theo}{Theorem}
\mynewtheo{rem}{Remark}
\mynewtheo{defi}{Definition}
\mynewtheo{conj}{Conjecture}
\mynewtheo{corr}{Corollary}
\mynewtheo{prop}{Proposition}
\mynewtheo{question}{Question}
\mynewtheo{exam}{Example}

\begin{document}

\makeatletter

\newenvironment{myeqn*}[1]{\begingroup\def\@eqnnum{\reset@font\rm#1}%
\xdef\@tempk{\arabic{equation}}\begin{equation}\edef\@currentlabel{#1}}
{\end{equation}\endgroup\setcounter{equation}{\@tempk}\ignorespaces}

\newenvironment{myeqn}[1]{\begingroup\let\eq@num\@eqnnum
\def\@eqnnum{\bgroup\let\r@fn\normalcolor 
\def\normalcolor####1(####2){\r@fn####1#1}%
\eq@num\egroup}%
\xdef\@tempk{\arabic{equation}}\begin{equation}\edef\@currentlabel{#1}}
{\end{equation}\endgroup\setcounter{equation}{\@tempk}\ignorespaces}

\newcommand{\mybin}[2]{\text{$\Bigl(\begin{array}{@{}c@{}}#1\\#2%
\end{array}\Bigr)$}}
\newcommand{\mybinn}[2]{\text{$\biggl(\begin{array}{@{}c@{}}%
#1\\#2\end{array}\biggr)$}}

\def\overtwo#1{\mbox{\small$\mybin{#1}{2}$}}
\newcommand{\mybr}[2]{\text{$\Bigl\lfloor\mbox{%
\small$\displaystyle\frac{#1}{#2}$}\Bigr\rfloor$}}
\def\mybrtwo#1{\mbox{\mybr{#1}{2}}}

\def\myfrac#1#2{\raisebox{0.2em}{\small$#1$}\!/\!\raisebox{-0.2em}{\small$#2$}}

\def\myeqnlabel{\bgroup\@ifnextchar[{\@maketheeq}{\immediate
\stepcounter{equation}\@myeqnlabel}}

\def\@maketheeq[#1]{\def\theequation{#1}\@myeqnlabel}

\def\@myeqnlabel#1{%
{\edef\@currentlabel{\theequation}
\label{#1}\enspace\eqref{#1}}\egroup}

\def\rato#1{\hbox to #1{\rightarrowfill}}
\def\arrowname#1{{\enspace
\setbox7=\hbox{F}\setbox6=\hbox{%
\setbox0=\hbox{\footnotesize $#1$}\setbox1=\hbox{$\to$}%
\dimen@\wd0\advance\dimen@ by 0.66\wd1\relax
$\stackrel{\rato{\dimen@}}{\copy0}$}%
\ifdim\ht6>\ht7\dimen@\ht7\advance\dimen@ by -\ht6\else
\dimen@\z@\fi\raise\dimen@\box6\enspace}}

\def\epsfs#1#2{{\epsfxsize#1\relax\epsffile{#2.eps}}}

\def\@test#1#2#3#4{%
  \let\@tempa\go@
  \@tempdima#1\relax\@tempdimb#3\@tempdima\relax\@tempdima#4\unitxsize\relax
  \ifdim \@tempdimb>\z@\relax
    \ifdim \@tempdimb<#2%
      \def\@tempa{\@test{#1}{#2}}%
    \fi
  \fi
  \@tempa
}

\def\go@#1\@end{}
\newdimen\unitxsize
\newif\ifautoepsf\autoepsftrue

\unitxsize4cm\relax
\def\epsfsize#1#2{\epsfxsize\relax\ifautoepsf
  {\@test{#1}{#2}{0.1 }{4   }
		{0.2 }{3   }
		{0.3 }{2   }
		{0.4 }{1.7 }
		{0.5 }{1.5 }
		{0.6 }{1.4 }
		{0.7 }{1.3 }
		{0.8 }{1.2 }
		{0.9 }{1.1 }
		{1.1 }{1.  }
		{1.2 }{0.9 }
		{1.4 }{0.8 }
		{1.6 }{0.75}
		{2.  }{0.7 }
		{2.25}{0.6 }
		{3   }{0.55}
		{5   }{0.5 }
		{10  }{0.33}
		{-1  }{0.25}\@end
		\ea}\ea\epsfxsize\the\@tempdima\relax
		\fi
		}


\author{A. Stoimenow\footnotemark[1]\\[2mm]
\small Ludwig-Maximilians University Munich, Mathematics\\
\small Institute, Theresienstra\ss e 39, 80333 M\"unchen, Germany,\\
\small e-mail: {\tt stoimeno@informatik.hu-berlin.de},\\
\small WWW: {\hbox{\tt http://www.informatik.hu-berlin.de/%
\raisebox{-0.8ex}{\tt\~{}}stoimeno}}
}

{\def\thefootnote{\fnsymbol{footnote}}
\footnotetext[1]{Supported by a DFG postdoc grant.}
}

\title{\large\bf \uppercase{Some minimal degree Vassiliev invariants}\\[2mm]
\uppercase{not realizable by the Homfly and Kauffman polynomial}\\[4mm]
{\it\small This is a preprint. I would be grateful
for any comments and corrections!}}

\date{\large Current version: \today\ \ \ First version:
\makedate{29}{6}{1998}}

\maketitle

\makeatletter

\def\chrd#1#2{\picline{1 #1 polar}{1 #2 polar}}
\def\arrow#1#2{\picvecline{1 #1 polar}{1 #2 polar}}

\def\labch#1#2#3{\chrd{#1}{#2}\picputtext{1.3 #2 polar}{$#3$}}
\def\labar#1#2#3{\arrow{#1}{#2}\picputtext{1.3 #2 polar}{$#3$}}
\def\labbr#1#2#3{\arrow{#1}{#2}\picputtext{1.3 #1 polar}{$#3$}}

\def\labline#1#2#3#4{\picvecline{#1}{#2}\pictranslate{#2}{
  \picputtext{#3}{$#4$}}}
\def\lablineb#1#2#3#4{\picvecline{#1}{#2}\pictranslate{#1}{
  \picputtext{#3}{$#4$}}}

\def\pt#1{{\picfillgraycol{0}\picfilledcircle{#1}{0.06}{}}}
\def\labpt#1#2#3{\pictranslate{#1}{\pt{0 0}\picputtext{#2}{$#3$}}}

\def\CD#1{{\let\@nomath\@gobble\small\diag{6mm}{2}{2}{
  \pictranslate{1 1}{
    \piccircle{0 0}{1}{}
    #1
}}}}

\def\GD#1{{\let\@nomath\@gobble\scriptsize\diag{6mm}{3.0}{3.0}{
  \pictranslate{1.5 1.5}{
    \piccircle{0 0}{1}{}
    #1
}}}}

\let\point\pt
\let\ay\asymp
\let\pa\partial
\let\al\alpha
\let\be\beta
\let\Gm\Gamma
\let\gm\gamma
\let\de\delta
\let\dl\delta
\let\eps\epsilon
\let\lm\lambda
\let\Lm\Lambda
\let\sg\sigma
\let\vp\varphi
\let\om\omega

\let\sm\setminus
\def\tl{\raisebox{-0.8 ex}{\tt\~{}}}

\def\ncap{\not\mathrel{\cap}}
\def\cf{\text{\rm cf}\,}
\def\lra{\longrightarrow}
\def\so{\Rightarrow}
\def\So{\Longrightarrow}
\let\ds\displaystyle
\def\bt{\bar t_2}

\let\reference\ref

\long\def\@makecaption#1#2{%
   \vskip 10pt
   {\let\label\@gobble
   \let\ignorespaces\@empty
   \xdef\@tempt{#2}%
   }%
   \ea\@ifempty\ea{\@tempt}{%
   \setbox\@tempboxa\hbox{%
      \fignr#1#2}%
      }{%
   \setbox\@tempboxa\hbox{%
      {\fignr#1:}\capt\ #2}%
      }%
   \ifdim \wd\@tempboxa >\captionwidth {%
      \rightskip=\@captionmargin\leftskip=\@captionmargin
      \unhbox\@tempboxa\par}%
   \else
      \hbox to\captionwidth{\hfil\box\@tempboxa\hfil}%
   \fi}%
\def\fignr{\small\sffamily\bfseries}%
\def\capt{\small\sffamily}%

\newdimen\@captionmargin\@captionmargin2cm\relax
\newdimen\captionwidth\captionwidth\hsize\relax

\def\eqref#1{(\protect\ref{#1})}

\def\proof{\@ifnextchar[{\@proof}{\@proof[\unskip]}}
\def\@proof[#1]{\noindent{\bf Proof #1.}\enspace}

\def\hint{\noindent Hint: }
\def\problem{\noindent{\bf Problem.} }

\def\@mt#1{\ifmmode#1\else$#1$\fi}
\def\qed{\hfill\@mt{\Box}}
\def\qqed{\hfill\@mt{\Box\enspace\Box}}

\def\cU{{\cal U}}
\def\cC{{\cal C}}
\def\cP{{\cal P}}
\def\tP{{\tilde P}}
\def\tZ{{\tilde Z}}
\def\fg{{\frak g}}
\def\tr{\text{tr}}
\def\cZ{{\cal Z}}
\def\cD{{\cal D}}
\def\bR{{\Bbb R}}
\def\cE{{\cal E}}
\def\bZ{{\Bbb Z}}
\def\bN{{\Bbb N}}

\def\br#1{\left\lfloor#1\right\rfloor}
\def\BR#1{\left\lceil#1\right\rceil}

\def\abstractname{}

\@addtoreset {footnote}{page}

\renewcommand{\section}{%
   \@startsection
         {section}{1}{\z@}{-1.5ex \@plus -1ex \@minus -.2ex}%
               {1ex \@plus.2ex}{\large\bf}%
}
\renewcommand{\@seccntformat}[1]{\csname the#1\endcsname .
\quad}

\def\bC{{\Bbb C}}
\def\bP{{\Bbb P}}

{\let\@noitemerr\relax
\vskip-2.7em\kern0pt\begin{abstract}
\noindent{\bf Abstract.}\enspace
We collect some examples showing that
some Vassiliev invariants are not obtainable from the HOMFLY and
Kauffman polynomials in the real sense, namely, that they distinguish
knots not distinguishable by the HOMFLY and/or Kauffman polynomial.
\\[1mm]
\end{abstract}
}

\section{Introduction}

Briefly after the introduction of Vassiliev knot invariants \cite{Vassiliev},
it was discovered \cite{BirmanLin}, that a lot of such invariants
can be ontained from knot polynomials \cite{HOMFLY,Jones,Kauffman},
basically by taking some coefficient in a version of the polynomial
where some variable was replaced by an exponential power series. In a
subsequent series of papers, several authors (see \cite{McDanielRong}
and \em{loc. cit.} and \cite{Dasbach}) independently computed the
number of Vassiliev invariants of given degree arising in this way
from the HOMFLY and Kauffman polynomials. This, however, is \em{a
priori}, only a lower bound for the number of Vassiliev invariants
of given degree ontainable from the knot polynomials as it is not
clear that no more can be generated in some other way. For example,
one could take some long polynomial expression of higher degree
Vassiliev invariants coming from the polynomial and little seems to be
known about how to control from below the degree of the resulting
Vassiliev invariant. Of course, one of the intrinsics of Vassiliev
theory (see \cite{BarNatanVI}) is that Vassiliev invariants form a
symmetric (i.e., commutative polynomial) algebra over a (primitive)
subspace of them, but trying to prove general degree bounds by looking
for the expression of a given Vassiliev invariant as a polynomial in
some fixed primitive basis is practically impossible, because
already the question of the number of primitive Vassiliev invariants
of given degree is far from being solved and one of the most
intriguing in the whole area (see \cite{ChmDuz,Dasbach2,Kneissler,%
Zagier}). See, however, \cite[\S 3]{restr}.
 
Given this subtlety, in this note we collect some examples showing that
some Vassiliev invariants are not ontainable from the HOMFLY and
Kauffman polynomials in the real sense, namely, that they distinguish
knots not distinguishable by the HOMFLY and/or Kauffman polynomial.
We hope these examples to be interesting at least because
(to the best of my knowledge) they never appeared explicitely
before in this context. All these examples are in a way related to
problems of detecting chirality with the knot polynomials.

\noindent{\bf Acknowledgement.} I would wish to thank to W.\ B.\ R.\
Lickorish, H.~Morton, D. Zagier, L.~Kauffman, T.~Kanenobu and
P.~Johnson for their helpful remarks.

\section{A degree 5 Vassiliev invariant not obtainable from the HOMFLY
polynomial}

It is known that the HOMFLY polynomial can be used to obtain all
Vassiliev invariants of degree at most 4 and two of the three 
primitive Vassiliev invariants of degree 5.

\begin{figure}[htb]
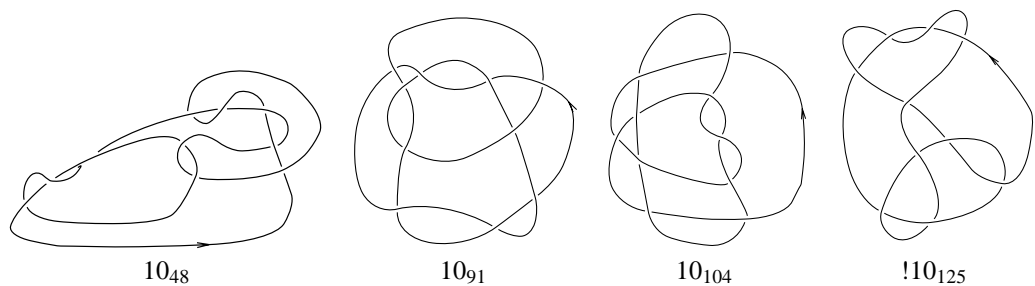

{
\[
\unitxsize3cm\relax
\begin{array}{*4c}
\epsfs{2cm}{k-10-48} &
\epsfs{2cm}{k-10-91} &
\epsfs{2cm}{k-10-104} &
\epsfs{2cm}{k-10-mi125} \\
10_{48} & 10_{91} & 10_{104} & !10_{125}
\end{array}
\]
}
\caption{The 4 knots with 10 crossings with 
self-conjuagate HOMFLY polynomial but non-self-conjuagate Kauffman 
polynomial.\label{4.10}}
\end{figure}

There are 4 knots with 10 crossings with
self-conjuagate HOMFLY polynomial but non-self-conjuagate Kauffman
polynomial, see figure \reference{4.10}. Each one of them together
with its obverse provide an example of two knots not distinguishable
by the HOMFLY polynomial but distinguishable by a degree 5 Vassiliev
invariant, showing that the third primitive Vassiliev invariant
of degree 5 is not obtainable from the HOMFLY polynomial.

\section{A degree 6 Vassiliev invariant not obtainable from the Kauffman
polynomial}

P. Johnson \cite{Johnson} used the tables \cite{polynomial} to obtain all
Vassiliev invariants of degree at most 5 from the Kauffman polynomial
and also 4 of the 5 primitive Vassiliev invariants of degree 6, but not the
fifth.

\begin{figure}[htb]
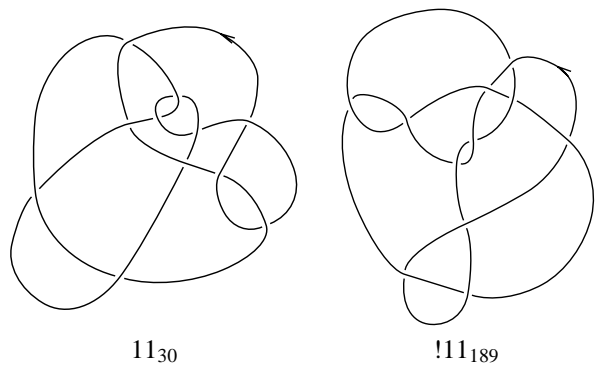

{
\[
\begin{array}{*2c}
\epsfs{6cm}{k-11-30} &
\epsfs{6cm}{k-11-mi189} \\
11_{30} & !11_{189}
\end{array}
\]
}
\caption{The pair of 11 crossing knots with
equal Kauffman polynomials but different Conway polynomials.
\label{2.11}}
\end{figure}

The example that the Kauffman polynomial does not contain a
degree 6 Vassiliev invariant is given by the two 11 crossing alternating
knots with equal Kauffman polynomials but different Conway polynomials,
see figure \reference{2.11}. In Thistlethwaite's notation \cite{Thi}
the knots are $11_{30}$ and $!11_{189}$. They Conway polynomials
are $\nabla(11_{30})=-2x^6 +x^4 -x^2 +1$ and $\nabla(!11_{189})=x^8 +2x^6 +x^4 -x^2 +1$, which differ at the coefficient of $x^6$, which is a
degree 6 Vassiliev invariant \cite{BarNatanVI}. This shows that
the Kauffman polynomial misses a degree 6 Vassiliev invariant,
which is contained in the Conway polynomial and hence also in the
HOMFLY polynomial. Thus, the HOMFLY and Kauffman polynomials
together exhaust all Vassiliev invariants up to degree 6.

\begin{rem}
The pair of knots $11_{30}$ and $!11_{189}$ was already known to
Lickorish at the 1986 Santa Cruz conference \cite[remark on
p.~472 top]{Kanenobu},
probably from the 11 crossing knot (polynomial) tabulation,
but has been reconstructed in a more elucidative way by Kanenobu
\cite[theorem 5]{Kanenobu}. Kanenobu observed that these knots can be used to
provide a counterexample to Lou Kauffman's hope \cite[p.~428]{Kauffman}
that the Kauffman polynomial always detects chirality when the HOMFLY
polynomial does so. Consider the knot $11_{30}\#11_{189}$. Then it has
self-conjuagate Kauffman polynomial but non-self-conjuagate HOMFLY
polynomial. Note also, that whenever the signature of a knot is not divisible
by 4, it has negative determinant $\Delta(-1)$, which is clearly
exhibited by the HOMFLY polynomial. So the HOMFLY polynomial detects chirality
in all such cases even if it is self-conjuagate, as for $9_{42}$
(where the Kauffman polynomial is also self-conjuagate) and $10_{125}$,
although not being able to distinguish between the knot and its obverse.
Therefore, care should be taken that the formulations ``$K$ has
self-conjuagate polynomial'' and ``the polynomial does not detect
chirality of $K$'' are not quite equivalent!
\end{rem}

\section{A degree 7 Vassiliev invariant not obtainable from the Kauffman
and HOMFLY polynomial}

Soon after the discovery of the three knot polynomials
capable of detecting chirality of knots, the knot $9_{42}$
became prominent by hiding its chirality to all of them.
A less noted relative of $9_{42}$ is $10_{71}$, which has
a crossing more but is more interesting because it has
additionally zero signature (contrarily to $9_{42}$).

Morton and Short \cite{MortonShort} developed a program
allowing to exhibit $9_{42}$'s chirality by showing
non-symmetry of the HOMFLY polynomial of its untwisted 2-cable.
An attempt to apply the same procedure on $10_{71}$
meets some (resource requirement) difficulties, as $10_{71}$
has braid index 5 (see \cite[appendix]{Jones}). Hence,
\em{a priori}, one has only a 10-braid representation
of its untwisted 2-cable. This was revealed to de indeed
minimal by the Morton-Williams-Franks inequality
(see \cite{Jones}) and the result of the calculation
shown on figure \reference{10-71}.
For its interpretation, if not self-explaining, we refer to
\cite{MortonShort}.

\begin{figure}[htb]
{\scriptsize
\begin{verbatim}

braid            : -2-3-1-287986576-4-5-3-4-4-5-3-421326576-4-5-3-48798-2-3-1-26576-4-5-3-4-6-7-5-68798                                                                            

algebraic cross no:  0

        -9        -7        -5        -3        -1         1         3         5         7         9
    ___________________________________________________________                 

         1        -7        21       -39        53       -53        39       -21         7        -1     | -1
        10       -57       146      -245       330      -365       308      -181        64       -10     |  1
        50      -246       540      -759       919     -1059       990      -645       260       -50     |  3
       112      -566      1182     -1467      1506     -1723      1789     -1294       573      -112     |  5
       113      -703      1552     -1823      1576     -1742      2042     -1606       704      -113     |  7
        54      -468      1226     -1447      1069     -1135      1525     -1238       468       -54     |  9
        12      -165       562      -716       461      -474       730      -563       165       -12     | 11
         1       -29       144      -208       120      -121       209      -144        29        -1     | 13
                  -2        19       -32        17       -17        32       -19         2               | 15
                             1        -2         1        -1         2        -1                         | 17

\end{verbatim}
}
\caption{(A part of) the output of the program of Morton and Short on
the untwisted 2-cable of $10_{71}$.\label{10-71}}
\end{figure}

Now, non-symmetry of the HOMFLY polynomial can be shown in terms of
Vassiliev invariants by evaluating at $i=\sqrt{-1}$ derivations
of the polynomial given by the rows of the output, i.~e. the
coefficients of the $z$ powers (in the convention of \cite{MortonShort}).
Then the $a$-th $v$-derivation at $v=i$ of the coefficients of $z^b$ is a
Vassiliev invariant of degree at most $a+b$ (by arguments analogous
to \cite{BirmanLin}). The difference between the polynomial
and its conjuagate comes out as a Vassiliev invariant of degree $7$
(e.g., for $b=1$ and $a=6$), which is therefore not contained neither
in the HOMFLY nor in the Kauffman polynomial, but in the
HOMFLY polynomial of the untwisted 2-cable.

\begin{figure}[htb]
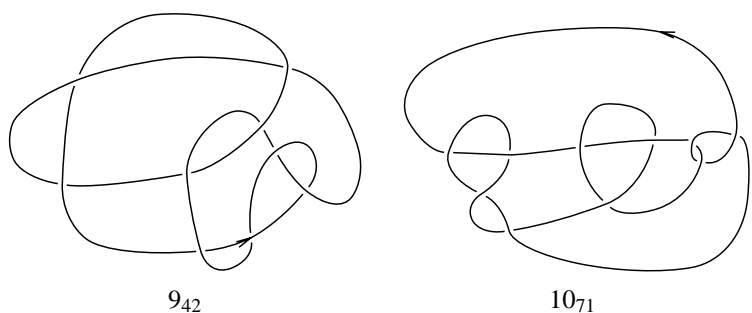

{
\[
\begin{array}{*2c}
\epsfs{6cm}{k-9-42} &
\epsfs{6cm}{k-10-71} \\
9_{42}  & 10_{71}
\end{array}
\]
}
\caption{The 2 chiral knots with at most 10 crossings with 
self-conjuagate HOMFLY and Kauffman polynomials.\label{2.10}}
\end{figure}

\begin{rem}
Przytycki conjectures that all Vassiliev invariants of
degree at most 10 are contained in the HOMFLY and Kauffman polynomials
and their untwisted 2-cables. A counterexample to this conjecture
would require a pair of knots not distinguishable by all of them. Such
pairs are mutants \cite{Conway,LicLip,Przytycki}. However, there is
no Vassiliev invariant known distinguishing mutants below degree 11
(see \cite{MorCro}). Are there more such pairs?
\end{rem}

\end{document}